\documentclass[12pt]{article}

\usepackage{amsthm,amsmath}
\usepackage{amssymb}
\usepackage{graphicx}
\usepackage{booktabs}
\def \qed {\hfill \hbox {\rule [-2pt]{3pt}{6pt}}}
\newtheorem{thm}{Theorem}
\newtheorem{theorem}[thm]{Theorem}

\newtheorem{lemma}[thm]{Lemma}
\newtheorem{remark}[thm]{Remark}

\newtheorem{definition}[thm]{Definition}

\def\Z{{\mathbf{Z}}}
\def\R{{\mathbf{R}}}
\def\N{{\mathbf{N}}}
\def\Z{{\mathbf{Z}}}

\begin{document}
\title{Algebraic Structure of Vector Fields \\in Financial Diffusion Models\\
and its Applications}
 
\author{
Yusuke MORIMOTO 
\thanks{
 The Bank of Tokyo-Mitsubishi UFJ, Ltd,
2-7-1 Marunouchi, Chiyoda-ku, Tokyo 100-8388, Japan, \
E-mail: yuusuke.morimoto@gmail.com} 
and Makiko SASADA
\thanks{
Graduate School of Mathematical Sciences, 
The University of Tokyo, 
3-8-1 Komaba, Meguro-ku, Tokyo 153-8914, Japan, \
E-mail: sasada@ms.u-tokyo.ac.jp}
}
\date{}

\maketitle

\begin{abstract}
High order discretization schemes of SDEs by using free Lie algebra valued random variables
are introduced by Kusuoka (\cite{K}, \cite{K2}), Lyons-Victoir (\cite{Cub}), 
Ninomiya-Victoir(\cite{NV}) and Ninomiya-Ninomiya(\cite{NN}). These schemes are called
KLNV methods.
They involve solving the flows of vector fields associated with SDEs and it is usually done by numerical methods.
The authors found a special Lie algebraic structure on the vector fields 
in the major financial diffusion models. Using this structure, 
we can solve the flows associated with vector fields analytically and efficiently.
Numerical examples show that our method saves the computation time drastically. 
\end{abstract}




\section{Introduction}
We consider $N$-dimensional Stratonovich stochastic differential equation 
\begin{equation*}
X(t,x) = x + \sum_{i=0}^{d} \int_{0}^{t} V_{i}(X(s,x))\circ dB^i(s)
\end{equation*}
where $N, d \geqq 1,$
$W_0$ $ = \{ w\in C([0,\infty );{\bf R}^d); \; w(0) = 0 \} ,$
${\cal F}$ be the Borel algebra over $W_0$ 
and $\mu$ be the Wiener measure on $(W_0,{\cal F}).$
Let $B^i:[0,\infty )\times W_0 \to {\bf R},$ $i=1,\ldots ,d,$ be given 
by $B^i(t,w) =w^i(t),$ $(t,w)\in [0,\infty )\times W_0.$
Then $ \{ (B^1(t), \ldots ,B^d(t) ; t \in [0,\infty ) \}$ 
is a $d$-dimensional Brownian motion.
Let $B^0(t) = t,$ $t \in [0,\infty ).$
Let $V_0,V_1,\ldots ,V_d$ 
$ \in C^{\infty}_b({\bf R}^N;{\bf R}^N).$
Here $C^{\infty}_b({\bf R}^N;{\bf R}^n)$ denotes 
the space of ${\bf R}^n$-valued smooth functions defined 
in ${\bf R}^N$ whose derivatives of any order are bounded.

Let us define a semigroup of linear operators $\{P_t\}_{t \in [0,\infty)}$ by
$$(P_t f)(x) = E[f(X(t, x))], \quad t \in [0,\infty), \  f \in C_b^{\infty} (\R^N).$$
 
 It is a crucial problem in various fields of applied science to approximate the expectation
$(P_Tf)(x)$ for a given time $T$ and function $f$ as fast and accurately as possible. The aim of 
this paper is to give a new method for the fast and accurate approximation of $(P_Tf)(x)$.

There are two approaches to this problem, the partial differential equation (PDE) based approach and
the simulation based approach. Notice that $u(t,x)=(P_t f)(x)$ satisfies the following equation
\begin{equation*} 
\frac{\partial }{\partial t}u(t,x) = Lu(t,x), \quad u(0, x) = f(x),
\end{equation*}
where the second order differential operator $L$ is given by
$ L = V_0 + \frac{1}{2}\sum_{i=1}^d V_{i}^2.$ Here we regard elements in $C^{\infty}_b({\bf R}^N;{\bf R}^N)$ 
as vector fields on ${\bf R}^N$ via
$$(V_{i} f )(x) = \sum_{j=1}^N V_{i}^j(x) \frac{\partial f}{\partial x_j}(x)
, \quad f\in C_b^{\infty}(\R^N).$$
In the PDE based approach, we solve this equation numerically. Whereas it works well when dimension $N$ is relatively small, it is prohibitively slow
in higher dimension. In such a case, the simulation based approach is the only practical method. 

The simulation based approach usually consists of two steps.
The first step is the time-discretization of stochastic differential equations using 
a set of random variables $\bar{X}^n(T,x)$, which approximate $X(T,x)$ in a certain sense as $n \to \infty$, whose
samples can be obtained by analytic or numerical methods.
If the stochastic differential equation is analytically solved, this 
step can be skipped. Otherwise, we apply a discretization scheme, such as Euler-Maruyama scheme which we will 
explain below, to get a random variables $\bar{X}(t,x)$ which is close to $X(t,x)$ for small $t>0$ in a certain sense.
Then, we construct a set of random variables $\{\bar{X}^n(t_k, x)\}_{k=0}^n$,
$(0 = t_0 < t_1 < \cdots < t_n = T )$ by repeating this approximation procedure $n$-times.
The second step is an approximation of $E[f(\bar{X}^n(T,x))]$ by Monte Carlo method (MC)
or Quasi Monte Carlo method (QMC). Both of them are essentially obtained  by averaging 
$M$ samples of $f(\bar{X}^n(T,x))$ denoted by $\{f(\bar{X}_m^n(T,x))\}_{m=1}^M.$
Samples are created at random for MC method and 
in a deterministic way for QMC method. The error is roughly estimated $O(M^{-1/2})$ for MC method
and $O(M^{-1})$ for QMC method.

In quantitative finance $X(t, x)$ represents the price of underlying assets 
and $E[f(X(T, x))]$ represents a price of derivative whose payoff function is $f$. 
For financial models, since the dimension of assets $N$ is often large, 
finding a fast and accurate method of the simulation based approach is very important.

In this paper, we consider the efficient scheme for the first step (discretization).
For a discretization scheme $\bar{X}(t,x)$, the linear operator $Q_t$ is defined by
$$(Q_tf)(x) = E[f(\bar{X}(t,x))].$$
Then the approximation of $(P_Tf)(x)$ with the approximation path $\{\bar{X}^n(t_k, x)\}_{k=0}^n$,
$(0 = t_0 < t_1 < \cdots < t_n = T )$ is described as
$$(Q_{t_n-t_{n-1}}\cdots Q_{t_2-t_1}Q_{t_1}f) (x).$$
We say that the discretization scheme is the weak approximation of order $r$ if there exists a constant $C > 0$ such that  
\begin{align*}
|(P_T f)(x) -  ((Q_{T/n})^n f)(x)| \leqq C n^{-r}
\end{align*}
for all $x \in \R^N$ and $f \in C_b^{\infty}(\R^N)$.

The most popular method of discretization schemes is the Euler-Maruyama scheme.
\begin{itemize}
\item Euler-Maruyama scheme
\begin{align*}
\bar{X}^{EM}(t,x)&=x+ \tilde{V}_0(x)t + \sum_{i=1}^d V_{i}(x)\sqrt{t}Z^{i},
\end{align*}
where $\tilde{V_0}^k = V_0^k+\frac{1}{2}\sum_{i=1}^d \sum_{j=1}^N V_{i}^j \frac{\partial}{\partial x_j}V_{i}^k$ for $k=1,2,\dots,N$
and $\{ Z^i \}_{i=1,\dots,d}$ is a family of independent $N(0,1)$ random variables.

\end{itemize}
It is known that the order of this scheme is $1$. 

Several higher order schemes have been studied by Kusuoka, Lyons, Ninomiya and Victoir where a free Lie algebra plays
an essential role. Note that the elements in $C^{\infty}_b({\bf R}^N;{\bf R}^N)$ 
are regarded as vector fields on ${\bf R}^N$. Then we can define the Lie bracket
$$[V_{\alpha}, V_{\beta}] = V_{\alpha}V_{\beta}-V_{\beta}V_{\alpha}, \quad 0 \leqq \alpha\leqq \beta \leqq d$$
where $[V_{\alpha}, V_{\beta}]$ is again an vector field on ${\bf R}^N$.
Let $\mathcal{L}$ denotes the Lie algebra generated by
$(\{V_0,V_1,\dots, V_d\}, [\cdot, \cdot])$. Using flows of $\mathcal{L}$ valued random variables,
Kusuoka (\cite{K}, \cite{K2}) and Lyons-Victoir (\cite{Cub}) introduced 
higher order discretization schemes. For any vector field $V \in C_b^{\infty}(\R^N; \R^N)$, the flow of $V$ is a
diffeomorphism $\exp(V) : \R^N \to \R^N$ given by $\exp(V) = u(1,x),$ where $u(t, x), t\geqq 0$
is the solution of the
following ODE
\begin{align*}
\begin{cases}
&\frac{du(t,x)}{dt} =V(u(t,x)), \quad t>0,\\
&u(0,x) = x.
\end{cases}
\end{align*}

Kusuoka showed in \cite{K,K2} that if there is a sequence of one parameter family of $\mathcal{L}$-valued random variables
$(\xi_1(t),\dots,\xi_{\ell}(t))_{t \geqq 0}$ satisfying some good condition with respect to $m \geqq 1$, a weak approximate operator
$Q^{(K)}_t$ is constructed with corresponding stochastic flows 
$\exp(\xi_1(t)), \dots, \exp(\xi_\ell(t))$ as
\begin{align*}
(Q^{(K)}_t f)(x) = E[f( \exp(\xi_\ell(t)) \circ \dots \circ \exp(\xi_1 (t)) ) (x))].
\end{align*}
More generally, if there is a set of sequences of one parameter family of $\mathcal{L}$-valued random variables
$(\xi_1^{i}(t),\dots,\xi_{\ell_i}^{i}(t))_{i=1,\dots,k}$ satisfying some good condition with respect to $m \geqq 1$, a weak approximate operator
$Q^{(K)}_t$ is constructed with corresponding stochastic flows as
\begin{align*}
(Q^{(K)}_t f)(x) = \sum_{i=1}^k c_i E[f(\exp(\xi^i_{\ell_i}(t)) \circ \dots \circ\exp(\xi_1^i(t))(x))]
\end{align*}
where $c_i$ are proper constants.
Under the assumption, the error arise from 1-step approximation is estimated by
$$|(P_tf)(x) - (Q^{(K)}_t f)(x)| \leqq C_{f, V}t^{\frac{m+1}{2}}$$
for $t \in (0,1]$ and the total error of the $n$-step approximation is estimated by 
$$|(P_Tf)(x) - ((Q^{(K)}_{T/n})^n f)(x)| \leqq C_{f, V}n^{-\frac{m-1}{2}}.$$ Namely, this is the weak approximation of order $\frac{m-1}{2}$.

Ninomiya and Victoir found a practical example for $m = 5$ in \cite{NV}. 
\begin{itemize}
\item Ninomiya-Victoir scheme 
\begin{align} \label{eq:nv}
&\bar{X}^{(NV)}(t,x)  \\
=& 
\begin{cases}
\exp(\frac{t}{2}V_0)\circ \exp(\sqrt{t}Z^1V_1) \circ \dots \nonumber \\
\qquad	 \dots \circ \exp(\sqrt{t}Z^dV_d)\circ
 \exp(\frac{t}{2}V_0)(x), \ \text{if} \ N = 1, \nonumber\\
\exp(\frac{t}{2}V_0)\circ \exp(\sqrt{t}Z^dV_d)\circ \dots \nonumber \\
\qquad \dots \circ \exp(\sqrt{t}Z^1V_1)\circ
\exp(\frac{t}{2}V_0)(x), \ \text{if} \ N = -1 \nonumber
\end{cases}
\end{align}
where $N$ is a Bernoulli random variable with the distribution $P(N=1)=P(N=-1)=\frac{1}{2}$ and $\{ Z^i \}_{i=1,\dots,d}$ is a family of independent $N(0,1)$ random variables where $N$ and $\{ Z^i \}_{i=1,\dots,d}$ are also independent.
\end{itemize}
The weak approximate operator
$Q^{(NV)}_t$ is given as 
\begin{align*}
& (Q^{(NV)}_t f)(x)\\
 & = \frac{1}{2} E[f(\exp(\frac{t}{2}V_0)\circ \exp(\sqrt{t}Z^1 V_1) \circ \dots \circ \exp(\sqrt{t}Z^d V_d)\circ \exp(\frac{t}{2}V_0)(x))]  \\
& + \frac{1}{2}E[f(\exp(\frac{t}{2}V_0)\circ \exp(\sqrt{t}Z^d V_d) \circ \dots \circ \exp(\sqrt{t}Z^1 V_1)\circ \exp(\frac{t}{2}V_0)(x))]. 
\end{align*}

Ninomiya and Ninomiya found another one parameter family of practical examples for $m= 5$ in \cite{NN}, Theorem 1.6.
\begin{itemize}
\item Ninomiya-Ninomiya scheme
\begin{align} \label{eq:nn}
& \bar{X}^{(NN)}(t,x)  \\
=& \exp\left(r t V_0+\sum_{i=1}^d S_{i}^1  \sqrt{t} V_i\right) \circ 
\exp\left((1-r) t V_0+\sum_{i=1}^d S_{i}^2 \sqrt{t}V_i\right)(x) \nonumber
\end{align}
where $S_i^1=r Z_i^1+\frac{1}{\sqrt{2}}Z_i^2$, $S_i^2=(1-r)Z_i^1 - \frac{1}{\sqrt{2}}Z_i^2$ and $(Z_i^{j})_{j=1,2,i=1,\dots,d}$ is a family of 
independent $N(0,1)$ random variables. $r \in \R$ is an arbitrary chosen fixed parameter.
\end{itemize}
The weak approximate operator
$Q^{(NN)}_t$ is given as 
\begin{align*}
& (Q^{(NN)}_t f)(x)\\
 & =  E[f(\exp\left(r V_0+\sum_{i=1}^d S_{i}^1  \sqrt{t} V_i\right) \circ 
\exp\left((1-r)tV_0+\sum_{i=1}^d S_{i}^2 \sqrt{t}V_i\right)(x))]  
\end{align*}
\begin{remark}
In \cite{NN}, a family of Gaussian random variables $S_i^j$ are characterized by a parameter $u \geqq \frac{1}{2}$ satisfying 
$r=\frac{\mp \sqrt{2(2u-1)}}{2}$ as
$E[S_j^i S_{j'}^{i'}]=R_{jj'}\delta_{ii'}$ where
\begin{align*}
R_{11}=u, \quad R_{12}=-u\mp\frac{\sqrt{2(2u-1)}}{2} \quad R_{22}=1+u\pm \sqrt{2(2u-1)}.
\end{align*}
We find that such a family of Gaussian random variables is constructed by $S_i^1=r Z_i^1+\frac{1}{\sqrt{2}}Z_i^2$ and $S_i^2=(1-r)Z_i^1 - \frac{1}{\sqrt{2}}Z_i^2$ where $(Z_i^{j})_{j=1,2,i=1,\dots,d}$ is a family of 
independent $N(0,1)$ random variables.
\end{remark}

Lyons and Victoir also introduced high order scheme by using a free Lie algebra
called Cubature on Wiener space (\cite{Cub}).
These methods are called KLNV method (Kusuoka, Lyons, Ninomiya and Victoir).

To compute an approximate value $(Q^{(K)}_t f)(x)$ such as $(Q^{(NV)}_t f)(x)$ and $(Q^{(NN)}_t f)(x)$, 
we need to compose the solution flows of random vector fields.
The composition is usually relying on the numerical ODE solvers such as
high order Runge-Kutta method.
For some lucky pair of vector fields and a scheme, ODE flows can be solved
analytically.
In such a case, the computation speed of the approximate value improves much faster.

In this paper, we show that if $\mathcal{L}$ has a special Lie algebraic structure (we call this condition Witt condition), then for any sequence of one parameter family of $\mathcal{L}$-valued random variables $(\xi_1(t),\dots,\xi_{\ell}(t))_{t \geqq 0}$ appearing in KLNV method and any order $m$, there exists a sequence of one parameter family of random variables $(p_0(t),p_1(t),\dots,p_k(t))_{t \geqq 0}$ such that
\begin{align*}
| & E[f(\exp(\xi_{\ell}(t)) \circ \dots \circ \exp(\xi_1(t))(x))] \\
 & -E[f(\exp(p_k(t) W_k) \circ \dots  \circ \exp(p_1(t) W_1) \circ\exp(p_0(t) W_0) (x))]|
 \leqq C_{f,m,W} t^{\frac{m+1}{2}}
\end{align*}
where the set $\{W_n, n\geqq 0\}$ is a special basis of $\mathcal{L}$.
Therefore, if there is an analytic solution of the flow of vector fields $W_n, n \geqq 0$, we can compute the 
approximate value $((Q_{T/n}^{(K)})^n f)(x)$ without using ODE solvers, and so have a faster approximation method.
We emphasize that our result enables us to a higher speed simulation for any scheme of KLNV method with any order, including Ninomiya-Victoir scheme and Ninomiya-Ninomiya scheme.

As an worthy application of the result, we show that our condition on the Lie algebraic structure and the existence of the analytic solution of vector fields $W_n$ are both satisfied for SABR model and Heston model, which are most important financial models. 
It has been known that if we apply Ninomiya-Victoir scheme to Heston model, the vector fields appearing in $Q^{(NV)}_t$ are analytically solved (\cite{NV}).
Bayer et al. \cite{Bay} suggested that, for SABR model if we rewrite the SDE using drifted Brownian motions and then apply Ninomiya-Victoir scheme,
the vector fields appearing in $Q^{(NV)}_t$ are also analytically solved. On the other hand, for both models,
to apply Ninomiya-Ninomiya scheme without using our new method, we need to use the numerical ODE solvers which increase the computation time.


We apply our result to Ninomiya-Ninomiya scheme in SABR model and do a numerical experiment. It shows that
our method is enough accurate and the computation time is highly saved.

\section{Notation and Main Result}
Now, we introduce a precise condition required for the family of vector fields $\{V_0,V_1,\dots,V_d\}$.
\begin{definition}(Witt condition)
We say that the family of vector fields $\{V_0,V_1,\dots,V_d\}$ satisfies the Witt condition if there exists a set of vector fields $\{W_n;n \geqq 0 \}$ satisfies 
\begin{itemize}
\item $V_i \in \text{span}\{W_n, n \geqq 0 \}$ for $i=0,1,\dots,d$ 
\item For any $n,m \geqq 0$, $[W_n, W_m]=\alpha_{nm}W_{n+m}$ with $\alpha_{nm} \in \R$.
\end{itemize}
\end{definition}

\begin{remark}
If $\{V_0,V_1,\dots,V_d\}$ satisfies Witt condition, then it is obvious that $\mathcal{L} \subset \text{span}\{W_n, n \geqq 0 \}$.
\end{remark}

\begin{remark}
The Witt algebra is the well-known Lie algebra with a basis $\{U_n; n \in \Z \}$ satisfying $[U_n,U_m]=(n-m)U_{n+m}$, so we name the above condition Witt condition.
\end{remark}

From now on, we assume $\{V_0,V_1,\dots,V_d\}$ satisfies Witt condition and so that $\mathcal{L} \subset \text{span}\{W_n, n \geqq 0 \}$. To state our main theorem, we introduce some notations. Let $\R[\lambda]=\{\sum_{k = 0}^n a_k \lambda^k; n \in \N_0, a_k \in \R \}$ be the polynomial ring of $\lambda$ where $\N_0=\{0,1,2,\dots\}$. Define the space $\R[\lambda]^*$
as a set of special $\R[\lambda]$ valued random variables as $\R[\lambda]^*=\{\sum_{k = 1}^n Z_k \lambda^k; n \in \N, Z_k \in L^{\infty,-} \}$. Here $L^{\infty,-}$ denotes the set of random variables having finite moments of all orders. Notice the sum of $k$ is taken from $1$. Define $\mathcal{L}^*$ as a set of special $\R[\lambda] \otimes \mathcal{L}$ valued random variables $\mathcal{L}^*=\{\sum_{k = 0}^n Q_k W_k; n \in \N_0, Q_k\in \R[\lambda]^*\}$. For $t \geqq 0$, define an operator $\Psi_t : \R[\lambda] \to \R$ as $\Psi_t(\sum_{k = 0}^n a_k \lambda^k)=\sum_{k = 0}^n a_k (\sqrt{t})^k$. Then the operator $\Psi_t$ is naturally extended to $\R[\lambda]^*$ and $\R[\lambda] \otimes \mathcal{L}$.

Next theorem is the main result of the paper.

\begin{theorem}\label{thm:main}
For any $m \geqq 1$ and for any sequence $(\xi_1,\xi_2,\dots,\xi_{\ell})$ in $\mathcal{L}^*$, there exists a number $k \in \N$ and a sequence $(P_0,P_1,\dots,P_k)$ in $\R[\lambda]^*$ such that
\begin{align*}
|E  [f( \exp(\Psi_t(\xi_{\ell})) \circ \exp(\Psi_t(\xi_2))) & \circ \dots \circ \exp(\Psi_t(\xi_1))(x))] \\
 -E[f( \exp(\Psi_t(P_k W_k))  \circ  \dots \circ \exp( & \Psi_t(P_1 W_1)) \circ \exp(\Psi_t(P_0 W_0))(x))]| \\
& \leqq C_{f,m,W} t^{\frac{m+1}{2}}
\end{align*}
for any $f \in C_b^{\infty}(\R^N)$ and $t \in (0,1]$
where $C_{f,m,W}$ is a constant depends only on $f, m$ and $W$.
\end{theorem}
We give a proof in Section \ref{section:decom}. The coefficients $P_0,P_1,\dots$ can be explicitly calculated from $(\xi_1,\xi_2,\dots,\xi_{\ell})$ with the recursive algorithm given in Section 4.
In Section 5, we apply our main result to the approximation of a price of European option for SABR model using Ninomiya-Victoir scheme and Ninomiya-Ninomiya scheme. The simulations result shows 
our method is properly accurate and reduces the computation time drastically compared to the existing methods.

\subsection{Application to financial models}

Here, we show that SABR model and Heston model satisfy Witt condition and the basis of their vector fields are analytically solvable. 
Therefore, we have a new approximation method for these models.

\subsubsection{SABR model}
The SABR model is given by 
\begin{align}\label{eq:sabr}
&dX_1(t,x ) = X_2(t,x)X_1(t,x)_{+}^{\beta} dB^1(t) , \\
&dX_2(t, x) = \nu X_2(t,x) (\rho dB^1(t) + \sqrt{1-\rho^2} dB^2(t))\nonumber
\end{align}
where $x_{+}=x \vee 0$.

Let $\{W^S_n, n\geqq 0\}$ be the set of vector fields on $\R^2$ defined as
\begin{align*}
& W^S_{n}=\frac{1}{1-\beta}{x_1}_+^{1-n(1-\beta)}x_2^{n}\frac{\partial}{\partial x_1}, \quad n\geqq 1, \\
& W^S_0=-x_2\frac{\partial}{\partial x_2}
\end{align*}
for fixed $\frac{1}{2} < \beta < 1$. It is easy to see that they satisfy 
\begin{align*}
[W^S_n,W^S_m]=(n-m)W^S_{n+m}.
\end{align*}

The vector fields in SABR model are given by
\begin{align}\label{eq:vecS}
&V_0=\frac{1}{2} \nu^2 W^S_0 + \frac{1}{2}(\beta - 1)  \nu  \rho W^S_1 +\frac{1}{2}  \beta  (\beta - 1)W^S_2\nonumber \\
&V_1 = - \nu  \rho W^S_0+(1 - \beta) W^S_1\\
&V_2=-\nu  \sqrt{1 - \rho^2}W^S_0 \nonumber
\end{align}
where $\nu, \beta$ and $\rho$ are model parameters (cf. \cite{Bay}).
Therefore, we can apply Theorem \ref{thm:main} to the model. 
Moreover, we have the explicit expression of the flow of $W^S_n$ as
\begin{align*}
\exp(t W^S_n)(x)= \left(\left(n x_2^n t+{x_1}_+^{n(1-\beta)}\right)_+^{\frac{1}{n(1-\beta )}},x_2 \right)
\end{align*}
for $n \geqq 1$ and $\exp(tW^S_0)(x)=(x_1,x_2 \exp(-t))$.

\subsubsection{Heston model}
The Heston model is given by
\begin{align*}
&dX_1(t,x ) = \mu X_1(t,x)dt +\sqrt{X_2(t,x)_+}X_1(t,x) dB^1(t) , \\
&dX_2(t, x) = \kappa(\theta - {X_2}(t,x)_+)dt +\xi \sqrt{X_2(t,x)_+} (\rho dB^1(t) + \sqrt{1-\rho^2} dB^2(t)).
\end{align*}
Let $\{M_n, n\geqq 0\}$ and $\{L_n, n\geqq 0\}$ be the set of vector fields on $\R^2$ defined as
\begin{align*}
& M_n=2x_1{x_2}_+^{-\frac{n}{2}+1}\frac{\partial}{\partial x_1}, \qquad L_n=2{x_2}_+^{-\frac{n}{2}+1}\frac{\partial}{\partial x_2}.
\end{align*}
It is easy to see that they satisfy the relation
\begin{align}\label{eq:lieheston}
& [M_n,M_m]= 0 \nonumber \\
& [M_n,L_m]=(n-2)M_{n+m} \\
& [L_n,L_m]=(n-m)L_{n+m} \nonumber.
\end{align}
Now, let $W^H_{2n}=L_{n}$ and $W^H_{2n+1}=M_{n}$ for $n \geqq 0$. Then, from (\ref{eq:lieheston}), we have
\begin{align*}
[W^H_n,W^H_m]=c_{nm}W^H_{n+m}
\end{align*}
where
\begin{align*}
c_{nm}= \begin{cases}
0 \quad n,m  \quad  \text{odd} \\
\frac{n-5}{2} \quad n \quad \text{odd},\quad  m \quad \text {even} \\
\frac{n-m}{2} \quad n,m \quad \text{even}.
\end{cases}
\end{align*}
The vector fields in Heston model are given by 
\begin{align*}
& V_0 =\frac{1}{2}(\mu-\frac{\xi\rho}{4})W^H_5 -\frac{1}{4}W^H_1+\frac{1}{2}(\kappa \theta -\frac{\xi^2}{4})W^H_4-\frac{\kappa}{2}W^H_0 \nonumber \\
& V_1=\frac{1}{2}W^H_3+\frac{\xi\rho}{2} W^H_2 \\
& V_2=\frac{\xi \sqrt{1-\rho^2}}{2}W^H_2 \nonumber
\end{align*}
where $\mu,\xi,\rho,\kappa$ and $\theta$ are models parameters (cf. \cite{Bay}). Therefore, we can apply Theorem \ref{thm:main} to the model.
Moreover, we have the explicit expression of the flow of $M_n$ and $L_n$ as
\begin{align*}
&\exp(tM_n)(x)=\left(x_1\exp(2{x_2}_+^{-\frac{n}{2}} t ), x_2 \right),\\
&\exp(tL_n)(x)=\left(x_1, \left(n t+{x_2}_+^{\frac{n}{2}}\right)_+^{\frac{2}{n}} \right).
\end{align*}

\begin{remark}
Though we assumed that all vector fields $V_0, \dots, V_d$ are smooth as common setting in the analysis 
of higher order weak approximation methods, the vector fields in these financial models are not in $C_b^{\infty}(\R^2;\R^2)$, 
On the other hand, practically, we apply KLNV method or Theorem \ref{thm:main} 
to them and numerical experiments serves well.
Though the condition under which we can justify approximation methods rigorously in mathematics is very strict, 
we guess the range that these methods work well practically is a bit larger. 
\end{remark}

\begin{remark}
If the initial state $X(0,x) \equiv (x_1,x_2)$ satisfies $x_1>0,x_2>0$ in SABR model or Heston model, then $X(t,x) \in \R_+^2 \ a.s.$ so we do not
need to extend the vector fields to the nonnegative area to consider the solution of the SDE. However, when we apply the KLNV method, since we consider random vector fields, we sometimes need to consider the solution of 
flows at negative time. For this reason, we extend the vector fields and the solution of the flow to the whole space $\R^2$. At critical points, such as $(0,x_2)$ for 
SABR model, the vector fields and their flows are not smooth, but we apply KLNV method and Theorem \ref{thm:main} to these models as the same reason as the last remark. 
\end{remark}

\section{Decomposition to the Flow of Base Vector Fields}\label{section:decom}

Let $A=\{w_0,w_1,\dots\}=\{w_i;i \in \N_0 \}$, be an alphabet, a set of letters, and $A^*$ be the set of words consisting of $A$ including 
the empty word which is denoted by $1$. For $w=w_{i_1}w_{i_2}\dots w_{i_k} \in A^*, i_j \geqq 0, j=1,\dots, k , k \geqq 1$, we define $\| w\|=i_1+i_2 + \dots i_k + \sharp \{ j \in \{1,2,\dots, k\} ; i_j=0\}$ and $\|1\|=0$.

Note that $\R[\lambda]$ is the polynomial ring of $\lambda$. 
Let $\tilde{d} :\R[\lambda] \to {\bf Z}_{\geqq 0}$ be defined for $p=\sum_{k= 0}^{n} a_k \lambda^k$ as
\begin{align*}
\tilde{d}(p)=
\begin{cases}
\min \{ k \geqq 0; a_k \neq 0\}, &p \neq 0,\\
\infty , &p =0.
\end{cases}
\end{align*} 
\\
For $M >0$, define
$$\R[\lambda]\langle A \rangle_M = \{ \sum_{w \in A^*} p_w w \ ; \ p_w \in \R[\lambda], \tilde{d}(p_w) M \geqq \|w \| , \sharp\{ w \in A^*; p_w \neq 0 \} < \infty \}$$
and 
$$\R[\lambda]\langle\langle A\rangle\rangle_M = \{ \sum_{w \in A^* } p_w w \ ; \ p_w \in \R[\lambda],  \tilde{d}(p_w)M \geqq \|w \| \}.$$ 
Since $\tilde{d}(p_w p_v)=\tilde{d}(p_w)+\tilde{d}(p_v)$ and $\|wv\|=\|w\|+\|v\|$ for any $w,v \in A^*$, it is easy to see that $\R[\lambda]\langle A \rangle_M$ and $\R[\lambda]\langle\langle A\rangle\rangle_M$ are rings. 

\begin{remark}
$M$ represents the biggest ratio between the order of $\sqrt{t}$ and the index of $W$ (we regard the index of $W_0$ as $1$) appearing in the vector fields obtained by KLNV method. For example, applying Ninomiya-Victoir scheme or Ninom to SABR model, we
have the following vector fields $tV_0=tc_{00}W_0+tc_{01}W_1+tc_{02}W_2$, $\sqrt{t}V_1=\sqrt{t}c_{10}W_0+\sqrt{t}c_{11}W_1$ and $\sqrt{t} V_1=\sqrt{t}c_{20}W_0$ with random coefficients where $c_{ij}$ are constants. Then, we take $M=1$ since
 $\lambda^2c_{00}w_0+ \lambda^2 c_{01}w_1+\lambda^2c_{02}w_2$, $\lambda c_{10}w_0+\lambda c_{11}w_1$ and $\lambda c_{20}w_0  \in \R[\lambda]\langle A \rangle_{M}$ with $M=1$. If we use Ninomiya-Victoir scheme or Ninom for Heston model, then we take $M=3$.
\end{remark}

\begin{lemma}
For any $u \in \R[\lambda]\langle A \rangle_M$, $\exp(u) \in \R[\lambda]\langle\langle A\rangle\rangle_M$.
\end{lemma}
{\it Proof.}
Since $\R[\lambda]\langle A \rangle_M$ is a ring, 
$u^n \in \R[\lambda]\langle A \rangle_M$ for any $u \in \R[\lambda]\langle A \rangle_M$ and $n \in \N$.
So it follows that $\exp(u) = \frac{1}{n !}\sum_{n=0}^{\infty}u^n \in\R[\lambda]\langle\langle A\rangle\rangle_M.$  
\qed

For $u= \sum_{w \in A^* } p_w w  \in \R[\lambda]\langle\langle A\rangle\rangle_M $, define
\begin{align*}
d(u)=\min \{\tilde{d}(p_w) ; w \in A^*\}.
\end{align*} \\
Define the operators $\tilde{j}_m : \R[\lambda] \to \R[\lambda]$ as
 $$\tilde{j}_m(\sum_{k = 0}^{n} a_k \lambda^k)=\sum_{k = 0}^{m \wedge n} a_k \lambda^k,$$ 
and $j_m : \R[\lambda]\langle\langle A\rangle\rangle_M \to \R[\lambda]\langle A \rangle_M$ as
$$j_m (\sum_{w \in A^* } p_w w )= \sum_{w \in A^*} \tilde{j}_m(p_w) w= \sum_{w \in A^*, \|w\| \leqq Mm} \tilde{j}_m(p_w) w.$$
It is obvious that if $d(u) \geqq m+1$, then $j_m(u)=0$. Also, $j_m (\exp(u))=j_m(\exp(j_m(u)))$ for any $u \in \R[\lambda]\langle A\rangle_M$.

Let $\Phi:\R[\lambda]\langle A \rangle_M \to \R[\lambda] \otimes \mathcal{DO}(\R^N)$ be a homomorphism give by 
\begin{equation*}
\Phi(1)=Id, \quad \Phi(w_{i_1}\dots w_{i_k})=W_{i_1} \dots W_{i_k}, \quad k \geqq 1, i_j \geqq 0, j=1,2,\dots,k
\end{equation*}
where $\mathcal{DO}(\R^N)$ is the set of smooth differential operators on $\R^N$.
Note that we assume $[W_n, W_m]=\alpha_{nm}W_{n+m}$ for any $n,m \geqq 0$.

For $k \in \N_0$, let $\mathcal{L}_1^{(k)} = \{ \sum_{i=k}^{\ell} p_i w_i \in \R[\lambda]\langle A \rangle_M \ ; p_i \in \R[\lambda], \ {\ell} \geqq k \}$. 
For $n \geqq 1$, define $\mathcal{L}_n^{(k)}$ as the linear space of homogeneous Lie polynomials of order $n$ for $\mathcal{L}_1^{(k)}$. 
Precisely 
\begin{align*}
 \mathcal{L}_n^{(k)} & = \{ L_n(u_1,u_2,\dots,u_n) \in \R[\lambda]\langle A \rangle_M; \\
 &  L_n \text{ is a Lie polynomial of order} \  n, u_i \in \mathcal{L}_1^{(k)}, i=1,2,\dots, n\}.
\end{align*}



\begin{lemma}\label{lem:phi}
For any $u \in \mathcal{L}_n^{(k)}$ and $n \geqq 2$, there exists $v \in \mathcal{L}_1^{(nk+1)}$ such that $\Phi(u)=\Phi(v)$.
\end{lemma}
{\it Proof.}
Because of the Lie algebraic structure of $\{W_i\}$, for any homogeneous Lie polynomial $L$ of order $n$, $L(W_{i_1},W_{i_2}, \dots,W_{i_n})=cW_{i_1+i_2+\dots+i_{n}}$ with some $c \in \R$ and if $i_1= i_2=\dots=i_n$, then $c=0$. Therefore, if $i_j \geqq k$ for $j=1,2,\dots,n$, $L(W_{i_1},W_{i_2}, \dots,W_{i_n})=cW_{p}$ with $p \geqq nk+1$ and $c \in \R$.
\qed

\begin{lemma}\label{lem:cut}
For any $u \in \mathcal{L}_n^{(k)}$, $j_m(u)=0$ if $n \geqq m+1$.
\end{lemma}
{\it Proof.}
Since $\| w_i \| \geqq 1$ for any $i \in \N_0$, if $u=\sum_{i=k}^{\ell} p_iw_i \in  \mathcal{L}_1^{(k)}$, then $\tilde{d}(p_i) M \geqq 1$ for all $i$ and therefore $d(u) \geqq 1$. Consequently, for $u \in \mathcal{L}_n^{(k)}$, $d(u) \geqq n$.
\qed

\begin{lemma}\label{lem:zassen}
For any $u=\sum_{i=k}^{\ell} p_i w_i \in \mathcal{L}_1^{(k)}$, there exists a sequence $u_2,u_3,\dots$ satisfying $u_i \in \mathcal{L}_i^{(k)}, i=2,3,\dots$ such that
\begin{align*}
& j_m (\exp(u)) \\
& =j_m( \exp(p_kw_k)\exp(p_{k+1}w_{k+1}) \dots \exp(p_{\ell}w_{\ell}) \exp(u_2) \exp(u_3) \dots \exp(u_{m})).
\end{align*}
\end{lemma}
{\it Proof.}
Note the Zassenhaus formula (cf. Theorem 2 in \cite{S})
\begin{align*}
&\exp(t(X_1+X_2+ \dots+X_K)) \\
 =&\exp(tX_1)\exp(tX_2)\cdots\exp(tX_K)\exp(t^2 C_2(X_1,X_2,\dots,X_K)) \cdots\\
& \qquad \qquad \cdots \exp(t^n C_n(X_1,X_2,\dots,X_K)) \cdots.
\end{align*} 
for any elements $(X_1,X_2,\dots,X_K)$ in a Lie algebra where $C_n$ are homogeneous Lie polynomials of order $n$. Applying this formula to our $u=\sum_{i=k}^l p_i w_i$ and notice Lemma \ref{lem:cut}, we have the required decomposition.
\qed

\begin{lemma}\label{lem:decompose}
For any $u \in \mathcal{L}_1^{(k)}$, there exists $N \in \N$ and a sequence $q_1,q_2,\dots, q_N \in \R[\lambda]$ and $i_j \in \{ k,k+1,\dots,Mm\}, \ j=1,2, \dots, N$ such that 
$$\Phi(j_m (\exp(u)))=\Phi(j_m( \exp(q_1w_{i_1})\exp(q_2 w_{i_2}) \dots \exp(q_N w_{i_N}) )).$$
\end{lemma}
{\it Proof.}
From Lemma \ref{lem:phi} and \ref{lem:zassen}, for any $u=\sum_{i=k}^{\ell} p_i w_i \in \mathcal{L}_1^{(k)}$, there exists a sequence $v_2,v_3,\dots$ satisfying $v_i \in \mathcal{L}_1^{(ki+1)}, i=2,3,\dots$ such that
\begin{align*}
& \Phi(j_m (\exp(u))) \\
&=\Phi(j_m( \exp(p_kw_k)\exp(p_{k+1}w_{k+1}) \dots \exp(p_{\ell}w_{\ell}) \exp(v_2) \exp(v_3) \dots \exp(v_m))).
\end{align*}
Then, we can apply this fact to $v_2$, $v_3 \dots, v_m$ repeatedly. After repeating the procedure finite times (at most $(m-1)^{Mm}$ times since $j_m(u)=0$ if $u \in \mathcal{L}_1^{(Mm+1)}$), we get the complete decomposition.
\qed

\begin{lemma}\label{lem:order}
For any $p_1,p_2 \in \R[\lambda]$ and $k \leqq l$, there exists $N \in \N$ and a sequence $q_1,q_2,\dots, q_N \in \R[\lambda]$ and $i_j \in \{ k+1,k+2,\dots,Mm\}, \ j=1,2, \dots, N$ such that
\begin{align*}
& \Phi(j_m (\exp(p_2 w_{\ell})\exp(p_1 w_k))) \\
& =\Phi(j_m( \exp(p_1 w_k)\exp(p_2 w_{\ell}) \exp(q_1w_{i_1})\exp(q_2 w_{i_2}) \dots \exp(q_N w_{i_N}) )).
\end{align*}
\end{lemma}
{\it Proof.}
It is a simple consequence of the Zassenhaus formula again. Namely, we have 
$$j_m (\exp(p_1 w_k + p_2 w_{\ell}))=j_m( \exp(p_2 w_{\ell})\exp(p_1 w_k) \exp(v_2) \exp(v_3) \dots \exp(v_m))$$ and
$$j_m (\exp(p_1 w_k + p_2 w_{\ell}))=j_m( \exp(p_1 w_k)\exp(p_2 w_{\ell}) \exp(z_2) \exp(z_3) \dots \exp(z_m))$$
for some $v_i,z_i \in \mathcal{L}_i^{(k)}, i=2,3,\dots,m$. Then, we have 
\begin{align*}
j_m( \exp(p_2 w_{\ell}) & \exp(p_1 w_k) )  =j_m( \exp(p_1 w_k)\exp(p_2 w_{\ell}) \exp(z_2) \exp(z_3) \dots \\
 & \dots \exp(z_m))\exp(-v_m) \exp(-v_{m-1}) \dots \exp(-v_2)).
\end{align*}
Then, applying Lemma \ref{lem:decompose}, we complete the proof. 
\qed

\begin{theorem}\label{thm:decomposition}
For any sequence $u_1,u_2,\dots,u_n \in \mathcal{L}_1^{(0)}$ and $m \geqq 1$, there exists a sequence $p_1,p_2,\dots p_{Mm} \in \R[\lambda]$ such that
\begin{align*}
\Phi(j_m(\exp(u_1)\exp(u_2) \dots \exp(u_n)) )=\Phi(j_m(\prod_{i=0}^{Mm}\exp(p_i w_i))) 
\end{align*}
where $\prod_{i=0}^{Mm}\exp(p_i w_i)= \exp(p_0w_0)\exp(p_1w_1) \dots \exp(p_{Mm}w_{Mm})$.
\end{theorem}
{\it Proof.}
Applying Lemma \ref{lem:decompose} to each $\exp(u_i)$ for $i=1,2,\dots,n$, we have a decomposition where the index of $w$ is not ordered. Then, applying Lemma \ref{lem:order} repeatedly, we can arrange the components in numerical order of the index of $w$.
\qed

For $t \geqq 0$, we define an operator $\Phi_t : \R[\lambda]\langle A \rangle_M \to \mathcal{DO}(\R^N)$ as $\Phi_t=\Psi_t \circ \Phi$.

\begin{lemma} \label{taylor}
For any $n \in \N$ and $u_1,u_2,\dots, u_n \in \mathcal{L}_{1}^{(0)}$, 
\begin{align*}
&|f(\exp (  \exp ( \Phi_t (j_m u_n))\circ \dots \circ \exp ( \Phi_t ( j_m u_2)) \circ \Phi_t (j_m u_1))(x)))  \nonumber \\
& - \sum_{0 \leqq  k_1+k_2+\dots k_n \leqq m}\frac{1}{k_1! \dots k_n !}\Phi_t((j_m u_1)^{k_1}\dots (j_m u_n)^{k_n}) f)(x)| \nonumber \\
&  \leqq  \sum_{k_1+k_2+\dots k_n = m+1}\frac{1}{k_1! \dots k_n !}\|\Phi_t(j_m u_1)^{k_1}\dots (j_m u_n)^{k_n}) f)\|_{\infty}, 
\end{align*}
for any $t \in (0,1].$
\end{lemma}
{\it Proof.}
It follows easily from Lemma 11 of \cite{K}.
\qed

\begin{definition}
For $m \geqq 1$, $t \geqq 0$ and a sequence $\xi=(\xi_1,\xi_2,\dots,\xi_n)$ in $\mathcal{L}^*$, let 
\begin{displaymath}
Q^{\xi,m}_{t}f=E[f(\exp ( \Phi_t ( j_m \xi_n )) \circ  \dots \circ \exp ( \Phi_t ( j_m \xi_2 ) ) \circ \exp(\Phi_t ( j_m \xi_1) ) (x))  ].
\end{displaymath}
\end{definition}

\begin{definition}
We say that sequences $\xi=(\xi_1,\xi_2,\dots,\xi_n)$ and $\zeta=(\zeta_1,\zeta_2,\dots \zeta_{\ell})$ in $\mathcal{L}^*$ are $m$-similar 
if
\begin{displaymath}
\Phi(j_m ( \exp ( \xi_1 ) \exp ( \xi_2 ) \dots \exp ( \xi_n ) ) ) =\Phi ( j_m  ( \exp (\zeta_1 ) \exp (\zeta_2 ) \dots \exp (\zeta_{\ell} ) )) \quad a.s.
\end{displaymath}
\end{definition}
\begin{theorem}\label{thm:main3}
For any $m \geqq 1$ and a sequence $\xi=(\xi_1,\xi_2,\dots,\xi_n)$ in $\mathcal{L}^*$, there exists a 
sequence $P_0,P_1,\dots, P_{Mm}$ in $\R[\lambda]^*$ such that $\xi$ and $\zeta=(P_0w_0,P_1w_1,\dots, P_{Mm} w_{Mm})$ are $m$-similar.
\end{theorem}
{\it Proof.}
First, note that since $\xi_i \in \mathcal{L^*}$, we can find $M >0$ such that $\xi_i \in \R[\lambda]\langle A \rangle_M \ a.s.$ for any $i=1,2,\dots,n$.
Then, we can apply Theorem \ref{thm:decomposition} to obtain $P_0,P_1,\dots,P_{Mm}$.
\qed

\begin{theorem}\label{thm:approximation}
Let $\xi=(\xi_1,\xi_2,\dots \xi_n)$ and $\zeta=(\zeta_1,\zeta_2,\dots \zeta_{\ell})$ be sequences in $\mathcal{L}^*$. 
If they are $m$-similar, then 
there is a constant $C>0$ such that
\begin{equation*}
\|Q^{\xi,m}_{t}f-Q^{\zeta,m}_{t}f \|_{\infty} 
\leqq C t^{(m+1)/2} \sup_{k=0,1,\dots,m+1, \ \alpha_1 ,\dots, \alpha_k \in \{0,1,\dots,Mm\} } \|W_{\alpha_1}\dots W_{\alpha_k}f\|_{\infty}
\end{equation*} 
for any $t \in (0,1]$ and for any $f \in C_b^{\infty}(\R^N)$.
\end{theorem}
{\it Proof.}
Notice that
\begin{align*}
&j_m ( \exp ( \xi_1 ) \exp ( \xi_2 ) \dots \exp ( \xi_n ) ) \\
&=\sum_{0 \leqq  k_1+k_2+\dots k_n \leqq m}\frac{1}{k_1! \dots k_n !} (j_m \xi_1)^{k_1} \dots (j_m \xi_n)^{k_n} +R 
\end{align*}
where $$R=  \sum_{0 \leqq  k_1+k_2+\dots k_n \leqq m}\frac{1}{k_1! \dots k_n !} j_m( (j_m \xi_1)^{k_1} \dots (j_m \xi_n)^{k_n}- (j_m \xi_1)^{k_1} \dots (j_m \xi_n)^{k_n}).$$
Therefore, by Lemma \ref{taylor}, we have
\begin{align*}
&|E[\Phi_t(j_m ( \exp ( \xi_1 ) \exp ( \xi_2 ) \dots  \exp ( \xi_n ))f(x)] -Q^{\xi,m}_{t}f(x)|\\
\leqq &E[\sum_{k_1+k_2+\dots k_n = m+1}\frac{1}{k_1! \dots k_n !}\|(\Phi_t((j_m \xi_1)^{k_1}\dots (j_m \xi_n)^{k_n})f)\|_{\infty}] +E[\| \Phi_t (R) f\|_{\infty}].
\end{align*}
Let   
$j_m \xi_i= \sum_{\ell=0}^{Mm} (\sum_{k=1}^m Z_{k,\ell}^{(i)}\lambda^k) w_{\ell}$ with 
random variables $ Z_{k,\ell}^{(i)}  \in L^{\infty,-}$, $k =1,\dots, m, \ \ell=0,\dots,Mm, \ i =1, \dots, n.$ Then, if $k_1+k_2+\dots+k_n=m+1$,
$$E[\|(\Phi_t((j_m \xi_1)^{k_1} \dots (j_m \xi_n)^{k_n})f)\|_{\infty}]$$
$$ = E[\| \left(\sum_{\ell=0}^{Mm} (\sum_{k=1}^m Z_{k,\ell}^{(1)}t^{k/2}) W_{\ell}\right)^{k_1}  \dots \left(\sum_{\ell=0}^{Mm} (\sum_{k=1}^m Z_{k,\ell}^{(n)}t^{k/2}) W_{\ell}\right)^{k_n})f)\|_{\infty}]$$
$$ \leqq  E[ \left(\sum_{\ell=0}^{Mm} (\sum_{k=1}^m |Z_{k,\ell}^{(1)}| t^{k/2}) \right)^{k_1} \dots \left(\sum_{\ell=0}^{Mm} (\sum_{k=1}^m |Z_{k,\ell}^{(n)}|t^{k/2})\right)^{k_n}) ] \|f\|_{W,m+1}$$ where
$$\|f\|_{W,m+1}=\sup_{ \alpha_1 ,\dots, \alpha_{m+1} \in \{0,1,\dots, Mm \}}  \|W_{\alpha_1}  \dots W_{\alpha_{m+1}} f\|_{\infty}.$$
Since $ Z_{k,\ell}^{(i)}  \in L^{\infty,-}$, there exists a constant $C>0$ such that $$E[ \left(\sum_{\ell=0}^{Mm} (\sum_{k=1}^m |Z_{k,\ell}^{(1)}| t^{k/2}) \right)^{k_1} \dots \left(\sum_{\ell=0}^{Mm} (\sum_{k=1}^m |Z_{k,\ell}^{(n)}|t^{k/2})\right)^{k_n} ) ] \leqq Ct^{m+1/2}$$ for $t \in (0,1]$.
Since $j_m R=0$, with similar estimates, there exists a constant $C>0$ such that
$$ E[\| \Phi_t (R) f\|_{\infty}] \leqq C t^{m+1/2}\|f\|_{W,\le m}$$
for $t \in (0,1]$ where
$$\|f\|_{W, \leqq m}=\sup_{k=0,1,\dots,m} \sup_{ \alpha_1 ,\dots, \alpha_{k} \in \{0,1,\dots, Mm \}}  \|W_{\alpha_1}  \dots W_{\alpha_{k}} f\|_{\infty}.$$
Combining these estimates, we have a constant $C_1 >0$ satisfying
\begin{align*}
& |E[\Phi_t(j_m (  \exp ( \xi_1 ) \exp ( \xi_2 )  \dots \exp ( \xi_n )  ) )f(x)] -Q^{\xi,m}_{t}f(x)| \\
& \leqq C_1 t^{(m+1)/2} (\|f\|_{W,m+1}+\|f\|_{W,\le m}).
\end{align*}

In the same way, there exists a constant $C_2 >0$ such that 
\begin{align*}
&|E[\Phi_t(\exp ( \zeta_1 ) \exp ( \zeta_2) \dots j_m (\exp ( \zeta_n )) )f(x)] -Q^{\zeta,m}_{t}f(x)|\\
\leqq &C_2 t^{(m+1)/2}  (\|f\|_{W,m+1}+\|f\|_{W,\le m}).
\end{align*}
Since $\xi$ and $\zeta$ are $m$-similar, we have the assertion.
\qed

\begin{theorem}\label{thm:main2}
Let $\xi=(\xi_1,\xi_2,\dots \xi_n)$ be a sequences in $\mathcal{L}^*$. 
Then, for any $m \geqq 1$, there exists a sequence $P_0,P_1,\dots, P_{Mm}$ in $\R[\lambda]^*$ such that 
\begin{equation*}
\|Q^{\xi,m}_{t}f-Q^{\zeta,m}_{t}f \|_{\infty} \leqq C_{f,m,W} t^{(m+1)/2} 
\end{equation*} 
for any $t \in (0,1]$ where $\zeta=(P_0w_0,P_1w_1,\dots, P_{Mm}w_{Mm})$.
\end{theorem}
{\it Proof.}
It follows from Theorem \ref{thm:main3} and \ref{thm:approximation} straightforwardly. 
\qed

Theorem \ref{thm:main} is a direct consequence of Theorem \ref{thm:main2}. 

\section{Algorithm of Decomposition}\label{sec:p}
In this section, we explain the explicit algorithm to obtain the decomposition to the base vector fields $\{W_n\}_{n \ge 0}$.
Applying Theorem \ref{thm:decomposition} with $n=1$ and $a_1,a_2,\dots,a_k \in \R$, for each $m \ge 1$, we can find $P_{i,m}(\lambda,{\bf a}) \in \R[\lambda]$ satisfying
$$j_m \exp(\lambda \sum_{i=0}^{k}a_iW_i)= \prod_{i=0}^{\infty} \exp(P_{i,m}(\lambda,{\bf a})W_i).$$
Moreover, from the way of the construction, we can find universal $P_{i}(\lambda, {\bf a}) \in \R[[\lambda]]$ satisfying $j_m P_{i}(\lambda, {\bf a})=P_{i,m}(\lambda,{\bf a})$ for
any $m \ge 1$. Here, $\R[[\lambda]]$ is the set of formal series of $\lambda$.
So, we calculate $P_i(\lambda,{\bf a}), i=0,\ldots,$ such that $P_i(0,{\bf a})=0$ formally satisfying
$$ \exp(\lambda \sum_{i=0}^{k}a_iW_i)= \prod_{i=0}^{\infty} \exp(P_{i}(\lambda,{\bf a})W_i) = \exp(P_{0}(\lambda,{\bf a})W_0)\exp(P_{1}(\lambda,{\bf a})W_1)\dots$$ where the precise meaning is that we showed in the last section.
For this purpose, we follow the clever method introduced in \cite{Ca} to compute the Zassenhaus formula efficiently.

Let
\begin{align}\label{eq:rn}
 R_n(\lambda, {\bf a}) =\exp(-P_n(\lambda,{\bf a})W_n) &  \exp(-P_{n-1}(\lambda,{\bf a})W_{n-1})\cdots  \\
& \cdots \exp(-P_0(\lambda,{\bf a})W_0)\exp(\lambda \sum_{i=0}^{k}a_iW_i) \nonumber
\end{align}
which should be equal to
$$\prod_{i=n+1}^{\infty} \exp(P_i(\lambda,{\bf a})W_i)=\exp(P_{n+1}(\lambda,{\bf a})W_{n+1})\exp(P_{n+2}(\lambda,{\bf a})W_{n+2})\cdots.$$
Let
$$F_n(\lambda,{\bf a})=(\frac{d}{d \lambda} R_n(\lambda, {\bf a}))R_n( \lambda,{\bf a})^{-1}.$$
Then, by (\ref{eq:rn})
$$F_n(\lambda,{\bf a})=-\frac{d}{d \lambda} P_n(\lambda, {\bf a})W_n+\exp(-P_n(\lambda, {\bf a})W_n)
(\frac{d}{d \lambda} R_{n-1}(\lambda, {\bf a}))R_n(\lambda, {\bf a})^{-1}.$$
Since the last term satisfies $$\exp(-P_n(\lambda, {\bf a})W_n)
(\frac{d}{d \lambda} R_{n-1}(\lambda, {\bf a}))R_n(\lambda, {\bf a})^{-1}$$
$$=\exp(-P_n(\lambda, {\bf a})W_n)
(\frac{d}{d \lambda} R_{n-1}(\lambda, {\bf a})) \left(\exp(-P_n(\lambda,{\bf a})W_n )R_{n-1}(\lambda, {\bf a})\right)^{-1}$$
$$=\exp(-P_n(\lambda, {\bf a})W_n)
(\frac{d}{d \lambda} R_{n-1}(\lambda, {\bf a}))R_{n-1}(\lambda, {\bf a})^{-1} (\exp(-P_n(\lambda,{\bf a})W_n )^{-1}$$
$$=\exp(-P_n(\lambda, {\bf a})W_n)
F_{n-1}(\lambda,{\bf a}) (\exp(-P_n(\lambda,{\bf a})W_n )^{-1},$$
we have
$$F_n(\lambda,{\bf a})=-\frac{d}{d \lambda} P_n(\lambda, {\bf a})W_n+ \exp(-P_n(\lambda, {\bf a})W_n)
F_{n-1}(\lambda,{\bf a}) (\exp(-P_n(\lambda,{\bf a})W_n )^{-1}$$
$$=\exp(-P_n(\lambda, {\bf a})W_n)
(F_{n-1}(\lambda,{\bf a}) -\frac{d}{d \lambda} P_n(\lambda, {\bf a})W_n) (\exp(P_n(\lambda,{\bf a})W_n ).$$
The last term can be rewritten as
$$\exp(ad_{-P_n(\lambda,{\bf a})W_n})(F_{n-1}(\lambda,{\bf a}) -\frac{d}{d \lambda} P_n(\lambda, {\bf a})W_n) $$
which is  well known formula
$$\text{e}^AB\text{e}^{-A}=\text{e}^{ad_A}B=\sum_{n\geqq 0} \frac{1}{n!}ad_A^n B,$$
with $ad_AB=[A,B], ad_A^jB=[A,ad_A^{j-1}B], ad_A^0B=B.$\\

On the other hand, since $R_n(\lambda, {\bf a}) =\prod_{i=n+1}^{\infty} \exp(P_i(\lambda,{\bf a})W_i)$,
\begin{align}\label{eq:fn}
& F_n(\lambda, {\bf a}) =\frac{d}{d\lambda}P_{n+1}(\lambda, {\bf a})W_{n+1}  \\
& +\exp(P_{n+1}(\lambda, {\bf a})W_{n+1}) \left(\frac{d}{d\lambda}\prod_{i=n+2}^{\infty}\exp(P_i(\lambda, {\bf a})W_i)  \right)
R_n(\lambda, {\bf a})^{-1} \nonumber \\
& =\frac{d}{d\lambda}P_{n+1}(\lambda, {\bf a})W_{n+1}+\sum_{i=n+2}^{\infty}
\exp(P_{n+1}(\lambda,{\bf a})W_{n+1})\cdots \exp(P_{i-1}(\lambda,{\bf a})W_{i-1}) \nonumber \\
&(\frac{d}{d\lambda}P_{i}(\lambda, {\bf a})W_{i}) \exp(-P_{i-1}(\lambda,{\bf a})W_{i-1})\cdots \exp(-P_{n+1}(\lambda,{\bf a})W_{n+1}) \nonumber \\
& =\frac{d}{d\lambda}P_{n+1}(\lambda, {\bf a})W_{n+1} +G_n(\lambda,{\bf a})\nonumber
\end{align}
where
$$G_n(\lambda,{\bf a})=\sum_{i=n+2}^{\infty}
\exp(Ad_{P_{n+1}(\lambda,{\bf a})W_{n+1}})\cdots \exp(Ad_{P_{i-1}(\lambda,{\bf a})W_{i-1}})(\frac{d}{d\lambda}P_{i}(\lambda, {\bf a})W_{i}).$$
Then we have the following relation.
\begin{align}
G_n(\lambda,{\bf a}) =F_n(\lambda,{\bf a}) - \frac{d}{d\lambda} P_{n+1}(\lambda, {\bf a})W_{n+1}, \label{G}\\
F_n(\lambda,{\bf a})=\exp(ad_{-P_n(\lambda,{\bf a})W_n})G_{n-1}(\lambda,{\bf a}). \label{F}
\end{align}
Using (\ref{G}) and (\ref{F}), we have an algorithm to obtain $P_n(\lambda,{\bf a})$.

For $F=\sum_j a_j W_j$, we define  $$\pi_{W_i}F=a_i.$$
From (\ref{eq:fn}) and Witt condition, we have $\pi_{W_{i}}F_n(\lambda,{\bf a})=0$ for $i \le n$ and
$$\pi_{W_{n+1}}F_n(\lambda,{\bf a})=\frac{d}{d \lambda}P_{n+1}(\lambda,{\bf a}).$$
Therefore, once $F_{n-1}(\lambda,{\bf a})$ is given, we can calculate $P_{n}(\lambda,{\bf a})$ using $\frac{d}{d \lambda}P_{n}(\lambda,{\bf a})=\pi_{W_{n}}F_{n-1}(\lambda,{\bf a})$ and $P_{n}(0,{\bf a})=0$.
Also, from  (\ref{G}), we can calculate $G_{n-1}(\lambda,{\bf a})$, and then using (\ref{F}), we can calculate $F_{n}(\lambda,{\bf a})$.
So, if we have $F_0(\lambda,{\bf a})$, by iterating this algorithm, we have $P_n(\lambda,{\bf a})$ for all $n \ge 0$:
$$F_0(\lambda,{\bf a}) \to P_1(\lambda,{\bf a}) \to G_0(\lambda,{\bf a}) \to F_1(\lambda,{\bf a}) \to \cdots \to P_n(\lambda,{\bf a}).$$

\subsection{Calculation of $F_0(\lambda,{\bf a})$}
First, notice that by the construction, $P_0(\lambda,{\bf a})=\lambda a_0$. Then,
$R_0=\exp(-\lambda a_0 W_0)\exp(\lambda \sum_{i=0}^{k}a_i W_i)$. Therefore,
\begin{align*}
& F_0(\lambda,{\bf a})=(\frac{d}{d \lambda}R_0(\lambda,{\bf a})) R_0(\lambda,{\bf a})^{-1}\\
&=-a_0 W_0+\exp(-\lambda a_0 W_0)(\sum_{i=0}^{k} a_i W_i)\exp(\lambda a_0 W_0) \\
& =-a_0 W_0+\exp(-\lambda a_0 W_0) a_0W_0\exp(\lambda a_0 W_0) +\exp(-\lambda a_0 W_0) (\sum_{i=1}^{k} a_i W_i) \exp(\lambda a_0 W_0) \\
& =\sum_{i=1}^{k} a_i \exp(Ad_{-\lambda a_0 W_0}) W_i
\end{align*}
On the other hand
\begin{align*}
 \exp(Ad_{-\lambda a_0 W_0}) W_i & = \sum_{n\ge 0}\frac{1}{n!}Ad_{-\lambda a_0 W_0}^nW_i \\
& =\sum_{n \ge 0}\frac{(-\lambda a_0)^n}{n!}   [W_0, [ \cdots,[ W_0, W_i]\cdots] \\
&=\sum_{n \ge 0}\frac{(-\lambda a_0)^n}{n!}   {\alpha_{0i}}^n W_i =\exp(-\alpha_{0i} \lambda a_0) W_i.
\end{align*}
So we have
$$F_0(\lambda,{\bf a})=\sum_{i=1}^{k} a_i \exp(-\alpha_{0i} \lambda a_0)  W_i.$$
Then we start the iteration. In particular, we have
$$P_1(\lambda, {\bf a})=-\frac{a_1}{\alpha_{0i}a_0}(\exp(-\alpha_{0i}\lambda a_0)-1).$$
\subsection{Example for SABR model}
We give the explicit formula of $P_i({\bf a}):=P_i(1,{\bf a}), i=0,1,2,3,4$ for SABR model where $k=\infty$ formally. When we apply Ninomiya-Victoir scheme or Ninomiya-Ninomiya scheme, we only need to take $k=2$, or equivalently $a_3=a_4=\dots=0$. Also, since $M=1$, to obtain an approximation of order $\frac{m+1}{2}$, we only need to calculate $P_i({\bf a}), i =0,1,\dots,m$ from Theorem \ref{thm:main2}. Moreover,
since all coefficients of the odd order terms of $\sqrt{t}$ have a symmetric distribution on $\R$ for Ninomiya-Victoir scheme or Ninomiya-Ninomiya scheme, to obtain the approximation of order $\frac{m+1}{2}$ with $m=5$,
we only need to use the terms $P_i({\bf a}), i =0,1,\dots,4$.
$$P_0({\bf a})=a_0, \quad P_1({\bf a})=\frac{a_1 }{ a_0 }(\exp(a_0)-1), \quad P_2({\bf a})=\frac{a_2}{2 a_0} (\exp(2 a_0)-1),$$
\begin{align*}
 P_3({\bf a})=&\frac{ \exp(a_0) -1} {6 a_0^2}(-a_1  a_2 - \exp(a_0)  a_1  a_2 + 2 \exp(2 a_0)  a_1  a_2 \\
                        &+ 2  a_0  a_3 + 2  \exp(a_0)  a_0  a_3 + 2  \exp(2  a_0)  a_0  a_3),
\end{align*}
\begin{align*}
 & P_4({\bf a})=\frac{1 }{12 a_0^3} (-1 + \exp(a_0) )  (a_1 ^2 a_2 + \exp(a_0) a_1^2 a_2 - 5 \exp(2  a_0) a_1^2 a_2 \\
                        &+ 3 \exp( 3 a_0) a_1^2 a_2- 2 a_0  a_1  a_3 - 2 \exp(a_0)  a_0  a_1  a_3 - 2 \exp( 2 a_0)  a_0  a_1  a_3 \\
                        &+ 6 \exp( 3 a_0)  a_0  a_1  a_3 + 3  a_0^2 a_4 + 3 \exp(a_0)  a_0^2  a_4+ 3 \exp(2  a_0) a_0^2 a_4 + 3\exp(3  a_0)  a_0^2 a_4).
\end{align*}

\section{Numerical results}

In this section we report the result of numerical experiment. We compare the following five schemes.
\begin{center}
\begin{tabular}{ll} \toprule
Lable & Method \\ \midrule
NN-analytic & Ninomiya-Ninomiya Scheme with analytic\\ &computation   proposed in this paper .\\
NV-analytic & Ninomiya-Victoir Scheme using analytic\\
& computation proposed in this paper\\
NN-Rk & Ninomiya-Ninomiya Scheme using Runge-Kutta method.\\
NV-Rk & Ninomiya-Victoir Scheme using Runge-Kutta method.\\ 
EM & Euler-Maruyama Scheme\\ \bottomrule
\end{tabular}
\end{center}

\subsection{Setting}
We used Quasi Monte Carlo method in particular  Sobol sequence for integration and 
we take $M=10^7$ for sampling number. 

We provide the results to the experiment with the SABR model given by (\ref{eq:sabr}).
For this experiment, parameters are chosen $\beta=0.9,$ $\nu = 1.0,$  $\rho = -0.7,$
and the initial value $x= (1.0, 0.3).$  
We choose a European call option with maturity $T=1.0$ and strike price $K=1.05.$
For simplicity we assume that the interest rate is zero.
This setting is the same as the experiment 4.1 in \cite{Bay}, and $0.09400046$ is used for the true result.
In the current experiment, we also use this value for true value.

We explain the detail of construction for NV-analytic scheme and NN-analytic scheme .
\subsection{NV-analytic scheme}
We define matrices $A= (a_{i,j})_{i=0,1,2,j=0,1,2}, B=(b_{i,j})_{i=0,1,2,j=0,1,2}$ and $C=(c_{i,j})_{i=0,1,2,j=0,1,2}$ as follows.
\begin{align*}
A=BC, \quad
B=
\begin{pmatrix}
\frac{t}{2}  & 0 & 0 \\
0 & \sqrt{t} Z^1 & 0 \\
0 & 0 & \sqrt{t} Z^2
\end{pmatrix},
\end{align*}
and
\begin{align}\label{eq:c}
C=
\begin{pmatrix}
\frac{1}{2} \nu^2 & \frac{1}{2}(\beta - 1)  \nu  \rho & \frac{1}{2}  \beta  (\beta - 1)\\
 - \nu  \rho & 1 - \beta & 0 \\
-\nu  \sqrt{1 - \rho^2} & 0 & 0 
\end{pmatrix},
\end{align}
where $Z^i, i = 1,2$ are independent $N(0, 1)$ random variables. 

Then, from the equation (\ref{eq:vecS}), we have
\begin{align*}
\begin{pmatrix}
\frac{t}{2}V_0\\
\sqrt{t}Z^1V_1\\
\sqrt{t}Z^2V_2
\end{pmatrix}
=A
\begin{pmatrix}
W_0\\
W_1\\
W_2
\end{pmatrix}.
\end{align*}
In the NV-analytic scheme, we decompose each flow $\exp(b_{ii} V_i)(x)$ in the equation (\ref{eq:nv}) as follows.
\begin{align*}
\exp(b_{ii} V_i)(x) = \exp(P_3 (a_{i0},a_{i1}, a_{i2})W_3 ) \circ \exp(P_2 (a_{i0},a_{i1}, a_{i2})W_2 ) 
\circ \\
\circ \exp(P_1 (a_{i0},a_{i1}, a_{i2})W_1 )\circ \exp(P_0 (a_{i0},a_{i1}, a_{i2})W_0 )(x) \nonumber
\end{align*}
where $P_i, i = 0, \dots , 3$ are the functions defined in Section \ref{sec:p}. 

\subsection{NN-analytic scheme}
Let $Y_0, Y_1$ be
\begin{align*}
&Y_0 = r t V_0+(r Z_1^1+\frac{1}{\sqrt{2}}Z_1^2)\sqrt{t} V_1+ (r Z_2^1 + \frac{1}{\sqrt{2}}Z_2^2) \sqrt{t}  V_2 ,\\
&Y_1 = (1-r) t V_0+((1-r) Z_1^1-\frac{1}{\sqrt{2}}Z_1^2)\sqrt{t}V_1+ ((1-r)Z_1^1 - \frac{1}{\sqrt{2}}Z_2^2) \sqrt{t}  V_2 ,
\end{align*}
where $Z_i^j, i, j = 1,2$ are independent $N(0, 1)$ random variables.
Applying NN-scheme (\ref{eq:nn})
\begin{align*}
\begin{pmatrix}
Y_0\\
Y_1
\end{pmatrix}
=A
\begin{pmatrix}
W_0\\
W_1\\
W_2
\end{pmatrix},
\end{align*}
where matrices $A= (a_{i,j})_{i=0,1,2,j=0,1}, B=(b_{i,j})_{i=0,1,j=0,1,2}$ and $C=(c_{i,j})_{i=0,1,2,j=0,1,2}$ are
\begin{align*}
A=BC, \quad
B=
\begin{pmatrix}
r t & (r Z_1^1+\frac{1}{\sqrt{2}}Z_1^2)\sqrt{t} & (r Z_2^1 + \frac{1}{\sqrt{2}}Z_2^2) \sqrt{t}  & \\
(1-r) t & ((1-r) Z_1^1-\frac{1}{\sqrt{2}}Z_1^2)\sqrt{t} & ((1-r)Z_1^1 - \frac{1}{\sqrt{2}}Z_2^2) \sqrt{t} 
\end{pmatrix}
\end{align*}
and
$C$ is defined by equation (\ref{eq:c}).

In the NV-analytic scheme, we decompose each flow $\exp(Y_i)(x)$ in equation(\ref{eq:nn})
as follows.
\begin{align*}
&\exp(Y_i)(x) = \exp(P_4 (a_{i0},a_{i1}, a_{i2})W_4) \circ \exp(P_3 (a_{i0},a_{i1}, a_{i2})W_3 ) \\
&\circ \exp(P_2 (a_{i0},a_{i1}, a_{i2})W_2) 
\circ \exp(P_1 (a_{i0},a_{i1}, a_{i2})W_1 )\circ \exp(P_0 (a_{i0},a_{i1}, a_{i2})W_0 )(x).
\end{align*}
where $P_i, i = 0, \dots , 4$ are the polynomials defined in Section \ref{sec:p}. 

\subsection{Result of experiment}
Figure 1 shows the convergence rates of five schemes.We see the both NN scheme and NV scheme
are actually second -order convergence and Euler scheme is first-order convergence.

Figure 2 shows the computation time of five schemes. 
The CPU used in this experiment is Intel(R) Core(TM) i7-46000 CPU@ 2.10GHz 2.7GH.
We can see the analytical scheme
save the computation time about $1/100$.
And NN analytic scheme and, NV analytic scheme are roughly the similar computation time.  

\begin{figure}
\begin{center}
\includegraphics[width=11cm, bb=0 0 673 643]{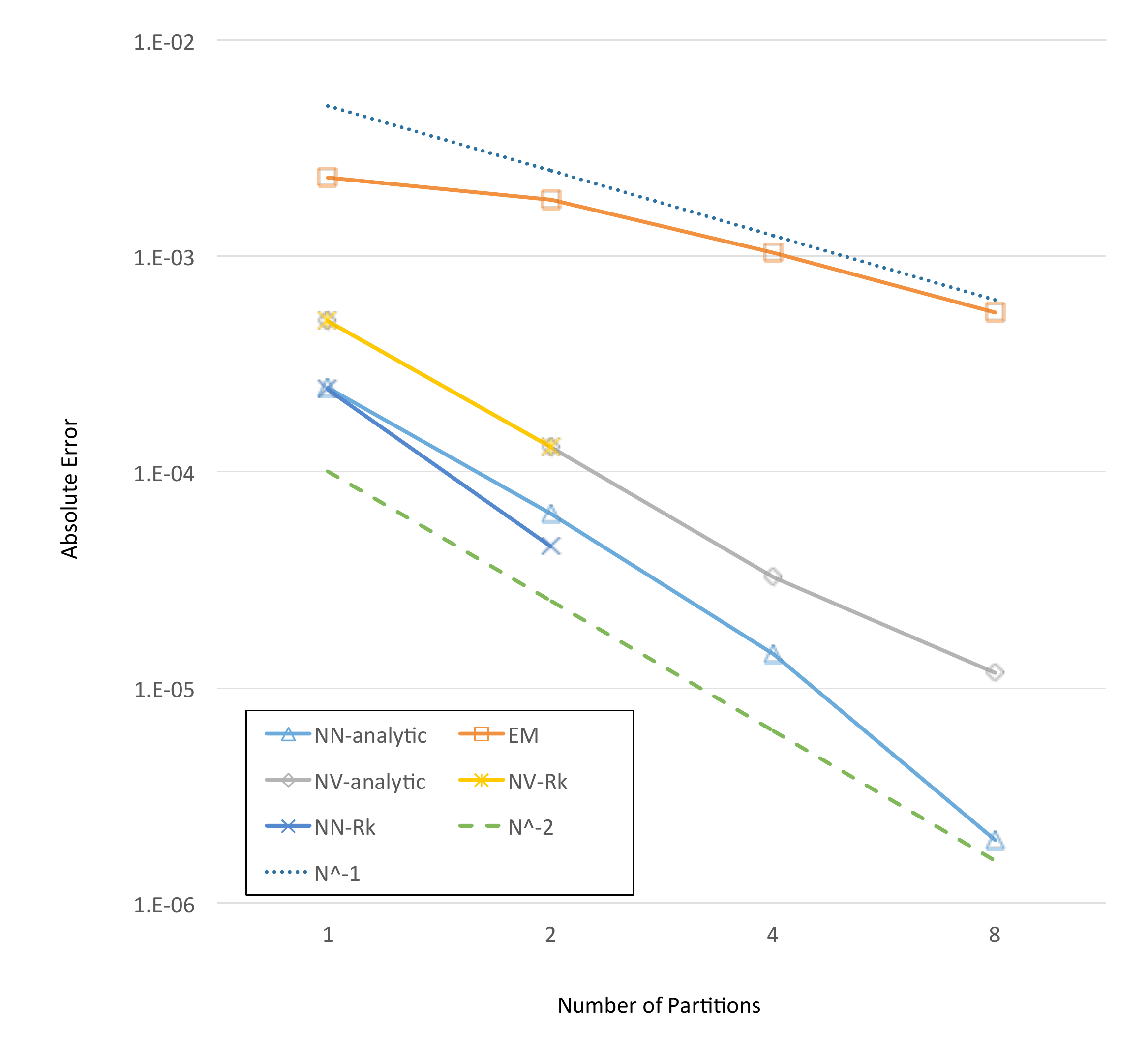}
\caption{Order of convergence of SABR model}
\end{center}
\bigskip
\begin{center}
\includegraphics[width=11cm, bb=0 0 668 537]{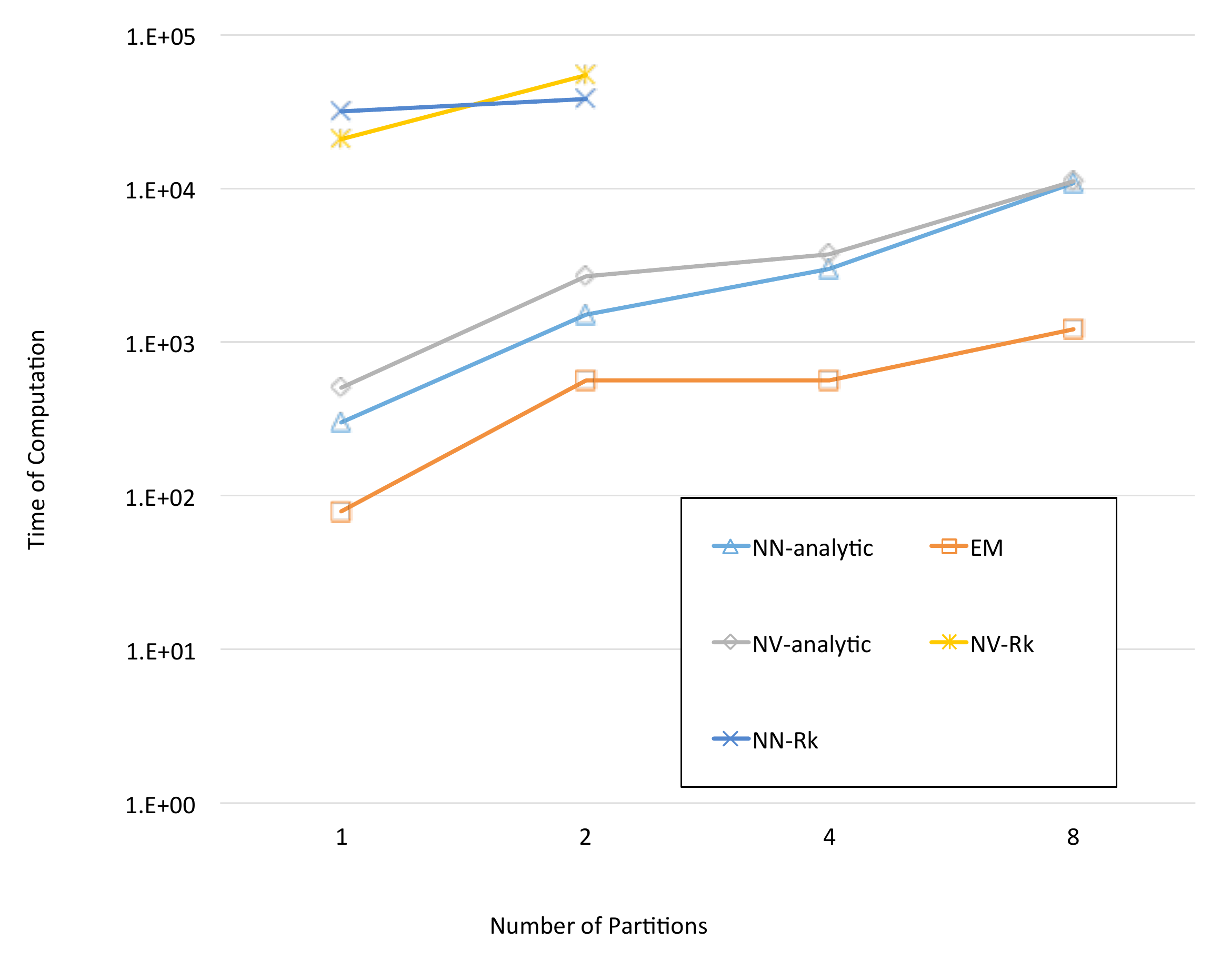}
\caption{Computation time of SABR model}
\end{center}
\end{figure}

\newpage

\end{document}